# ON ADMISSIBILITY
# FOR PARABOLIC EQUATIONS IN $\mathbb{R}^n$

### Martino Prizzi

Abstract. We consider the parabolic equation

(P)           $$u_t - \Delta u = F(x,u), \quad (t,x) \in \mathbb{R}_+ \times \mathbb{R}^n$$

and the corresponding semiflow $\pi$ in the phase space $H^1$. We give conditions on the nonlinearity $F(x,u)$, ensuring that all bounded sets of $H^1$ are $\pi$-admissibile in the sense of Rybakowski. If $F(x,u)$ is asymptotically linear, under appropriate non-resonance conditions, we use Conley's index theory to prove the existence of nontrivial equilibria of (P) and of heteroclinic trajectories joining some of these equilibria. The results obtained in this paper extend earlier results of Rybakowski concerning parabolic equations on *bounded* open subsets of $\mathbb{R}^n$.

## 1. Introduction

For $n \geq 3$, we consider the parabolic equation

$$(1.1) \qquad u_t - \Delta u = F(x,u), \quad (t,x) \in \mathbb{R}_+ \times \mathbb{R}^n.$$

The function $F \colon \mathbb{R}^n \times \mathbb{R} \to \mathbb{R}$ is assumed to be continuous; moreover, for every $x \in \mathbb{R}^n$, the function $F(x,\cdot)$ is assumed to be continuously differentiable. The assumption $n \geq 3$ is inessential and we make it only for notational convenience.

Associated with $-\Delta$, we consider the corresponding positive self-adjoint operator $A \colon D(A) \subset L^2(\mathbb{R}^n) \to L^2(\mathbb{R}^n)$, where $D(A) = H^2(\mathbb{R}^n)$, and for $u \in D(A)$, $Au = -\Delta u$ in the distributional sense. Recall that, for $s \in [0,1]$, $D((I+A)^s) = H^{2s}(\mathbb{R}^n)$. In particular, $H^1(\mathbb{R}^n) = D((I+A)^{1/2})$ and, for $u \in H^1(\mathbb{R}^n)$, $\|u\|_{H^1} = \|(I+A)^{1/2}u\|_{L^2}$. The operator $A$ generates an analytic semigroup of linear operators $e^{-tA}$, $t \geq 0$, for which the estimates

$$(1.2) \qquad \|(I+A)^s e^{-At}u\|_{L^2} \leq M\left(1 + \frac{1}{t^{(s-r)}}\right)\|(I+A)^r u\|_{L^2}, \quad t > 0,$$

hold. Here $M$ is a positive constant and $0 \leq r \leq s \leq 1$.

Let $C$ be a positive constant and set $\beta := 2/(n-2)$. If the conditions

$$(1.3) \qquad \begin{cases} F(\cdot, 0) \in L^2(\mathbb{R}^n) \\ |F'_u(x, u)| \leq C(1 + |u|^\beta), \quad (x, u) \in \mathbb{R}^n \times \mathbb{R}, \end{cases}$$

are satisfied, then we can define the Nemitski operator $\hat{F}(u)(x) := F(x, u(x))$ which, thanks to the Sobolev embedding $H^1 \hookrightarrow L^{2(\beta+1)}$, turns out to be a (non-linear) $C^1$ map of $H^1(\mathbb{R}^n)$ to $L^2(\mathbb{R}^n)$. Moreover, the estimates

$$(1.4) \qquad \|\hat{F}(u)\|_{L^2} \leq \tilde{C}(1 + \|u\|_{H^1}^{\beta+1})$$

and

$$(1.5) \qquad \|D\hat{F}(u)\|_{\mathcal{L}(H^1, L^2)} \leq \tilde{C}(1 + \|u\|_{H^1}^\beta)$$

hold (cf [8]); so, in particular, $\hat{F}$ is bounded and Lipschitz continuous on any bounded subset of $H^1(\mathbb{R}^n)$.

We can rewrite equation (1.1) in the abstract form

$$(1.6) \qquad \dot{u} + Au = \hat{F}(u).$$

Equation (1.6) generates a (local) semiflow $u\pi t$ in the phase space $H^1(\mathbb{R}^n)$ (see [6]). We recall also the variation-of-constant formula

$$(1.7) \qquad u(t) = e^{-At}u(0) + \int_0^t e^{-A(t-s)}\hat{F}(u(s))\, \mathrm{d}s.$$

The aim of this paper is to give conditions on $F$ ensuring that all bounded sets of $H^1(\mathbb{R}^n)$ are $\pi$-admissible (see Definition 2.1 below). In the case of a parabolic equation like (1.1) on a *bounded* open set $\Omega \subset \mathbb{R}^n$, the admissibility of all bounded subsets of $H^1(\Omega)$ is a direct consequence of the compactness of the Sobolev embedding $H^1(\Omega) \hookrightarrow L^2(\Omega)$. In $\mathbb{R}^n$ this property fails, and one has to introduce some restrictions on the non-linear term $F(x, u)$ (see (2.1) below). Roughly speaking, these restrictions mean that the nonlinearity $F(x, u)$ is dissipative for large $x$.

The concept of $\pi$-*admissible set* for a local semiflow $\pi$ in a metric space $X$ was introduced by Rybakowski in [11]. If $K$ is an isolated $\pi$-invariant set for which there exists a $\pi$-admissible isolating neighborhood $N$ (see [13] for the precise definitions of this and of the related concepts), then one can prove that there exists a special isolating neighborhood $\mathcal{B} \subset N$ of $K$, called an *isolating block*, which has the property that solutions of $\pi$ are 'transverse' to the boundary of $\mathcal{B}$. Letting $\mathcal{B}^-$ be the set of all points of $\partial\mathcal{B}$ the solutions through which leave $\mathcal{B}$ in positive time direction, and collapsing $\mathcal{B}^-$ to one point, we obtain the *pointed space* $\mathcal{B}/\mathcal{B}^-$ with the distinguished *base point* $p = [\mathcal{B}^-]$. It turns out that the homotopy type $h(\mathcal{B}/\mathcal{B}^-, [\mathcal{B}^-])$ of $(\mathcal{B}/\mathcal{B}^-, [\mathcal{B}^-])$ does not depend on the choice of $\mathcal{B}$. This means that $h(\mathcal{B}/\mathcal{B}^-, [\mathcal{B}^-])$ depends only on the pair $(\pi, K)$, and we write $h(\pi, K) := h(\mathcal{B}/\mathcal{B}^-, [\mathcal{B}^-])$. $h(\pi, K)$ is called the *homotopy index* of $(\pi, K)$. For two-sided flows on locally compact spaces, the homotopy index is due to Charles Conley ([3]) and therefore it is called the *Conley index*. In the case of a local semiflow $\pi$ in an arbitrary metric space $X$, the extended homotopy index theory was developed by Rybakowski (see [13]) and rests in an essential way on the notion of $\pi$-*admissibility*. The most important properties of the Conley index are the following: (a) if $h(\pi, K) \neq \underline{0}$, then $K \neq \emptyset$; (b) the homotopy index is invariant under continuation, in the sense that, roughly speaking, it remains constant along 'continuous' deformations of the pair $(\pi, K)$; (c) if $u$ is a hyperbolic equilibrium of Morse index $k$, then $h(\pi, \{u\}) = \Sigma^k$, where $\Sigma^k$ is the homotopy type of a $k$-dimensional pointed sphere.

Concerning equation (1.1), one main feature of the corresponding semiflow $\pi$ is its gradient-like structure. In fact, if $P(x, u) := \int_0^u F(x, s)\, \mathrm{d}s$, then

$$V(u) := \frac{1}{2} \int_{\mathbb{R}^n} |\nabla u(x)|^2 \,\mathrm{d}x - \int_{\mathbb{R}^n} P(x, u(x))\, \mathrm{d}x$$

is a Lyapunov functional for $\pi$, such that

$$\frac{d}{dt} V(u(t)) = -|\dot{u}(t)|^2$$

along any solution $u(t)$ of (1.1). Hence every nonempty compact $\pi$-invariant set $K$ contains at least one equilibrium of $\pi$, i.e. a solution of the elliptic equation

$$(1.8) \qquad\qquad -\Delta u = F(x, u), \quad x \in \mathbb{R}^n$$

(see e.g. [5], [7] or [16]). In this spirit, in Section 3 we shall consider equations which behave linearly at infinity. More precisely, we assume that $F(x, u)$ satisfies a condition like

$$\lim_{|u| \to \infty} \frac{F(x, u)}{u} = \alpha(x).$$

There is a wide literature on asymptotically linear elliptic equations in *bounded* domains, with or without resonance at infinity (see e.g. [2] and the references contained therein). We stress that various existence results have been obtained by a systematic use of the Conley index in [1] and later in [12], where (1.8) is considered in the context of the parabolic equation (1.1). Our second goal is to extend some results of [12] to equations on $\mathbb{R}^n$. As far as we know, very little has been done in the case of unbounded domains. We just quote the two recent papers [14], where the radial case is considered, and [4], where a very strong resonance at infinity is allowed, at the price of several severe restrictions on the behavior of $F(x, u)$ at finite. In both cases some 'mountain-pass' theorem is exploited. On the other hand, using Conley index techniques, we are somehow lead to consider only the non-resonance case, but a more general behavior of $F(x, u)$ at finite is allowed.

Like in the note [15] (where a saddle-point type theorem of Brezis and Nirenberg is used), we shall obtain an existence result for (1.8) under fairly general conditions on $F$. Besides, the topological information contained in the Conley index allow us to improve in a significant way the result of [15]: in fact, we are able to prove the existence of nontrivial solutions of (1.8) even when the techniques of [15] give no means to distinguish between trivial and nontrivial solutions. Finally, it is worth to mention that the *dynamical* approach (via Conley index theory) often gives as a by-product remarkable results on existence of heteroclinic trajectories joining some of the equilibria.

## 2. A CONDITION FOR ADMISSIBILITY

Whenever $\pi$ is a local semiflow in $X$, we write $x\pi t := \pi(t,x)$, $(t,x) \in \mathbb{R}_+ \times X$. We begin by recalling the following concept, introduced by Rybakowski in [11]:

**Definition 2.1.** *Let $X$ a metric space, let $N$ be a closed subset of $X$ and let $(\pi_j)_{j\in\mathbb{N}}$ be a sequence of local semiflows in $X$. Then $N$ is called $\{\pi_j\}$-admissible if the following holds:*

*if $(x_j)_{j\in\mathbb{N}}$ is a sequence in $X$ and $(t_j)_{j\in\mathbb{N}}$ is a sequence in $\mathbb{R}_+$ such that $t_j \to \infty$ as $j \to \infty$ and $x_j\pi_j[0,t_j] \subset N$ for all $j \in \mathbb{N}$, then the sequence of endpoints $(x_j\pi_jt_j)_{j\in\mathbb{N}}$ has a converging subsequence.*

*$N$ is called* strongly $\{\pi_j\}$-admissible *if $N$ is $\{\pi_j\}$-admissible and if $\pi_j$ does not explode in $N$ for every $j \in \mathbb{N}$. If $\pi_j = \pi$ for all $j$, we say that $N$ is $\pi$-admissible (resp.* strongly $\pi$-admissible*)*

Notice that, by [6, Th. 3.3.4], if $N \subset H^1(\mathbb{R}^n)$ is bounded then the semiflow $\pi$ generated by (1.1) does not explode in $N$.

In the case of a parabolic equation like (1.1) on a *bounded* open set $\Omega \subset \mathbb{R}^n$, the admissibility of all bounded subsets of $H^1(\Omega)$ is a direct consequence of the compactness of the Sobolev embedding $H^1(\Omega) \hookrightarrow L^2(\Omega)$. In $\mathbb{R}^n$ this property fails, and one has to introduce some restrictions on the non-linear term $F$. In this section we give conditions on $F$, ensuring that all bounded subsets of $H^1(\mathbb{R}^n)$ are $\pi$-admissible. Namely, we assume

$$(2.1) \qquad\qquad F(x,u)u \le -\nu|u|^2 + b(x)|u|^q + c(x)$$

where $\nu > 0$, $c \in L^1(\mathbb{R}^n)$, $2 \le q < 2n/(n-2)$ and $b \in L^p(\mathbb{R}^n)$, where $2n/[2n - q(n-2)] \le p < \infty$.

*Remark.* The results of this section still hold if in (2.1) one makes the alternative assumption that $2 \le q \le 2n/(n-2)$, $b \in L^\infty(\mathbb{R}^n)$ and $\lim_{k\to\infty}$ ess $\sup_{|x|\ge k} |b(x)| = 0$. In fact, one only needs to slightly modify the proof of Proposition 2.2 below. Roughly speaking, condition (2.1) means that the nonlinearity $F(x,u)$ is dissipative for large $x$.

Our first goal is to prove the following 'asymptotic localization' result, inspired by [17, Lemma 5]:

**Proposition 2.2.** *Assume $F(x,u)$ satisfies (1.3) and (2.1). Let $u\colon [0,T] \to H^1(\mathbb{R}^n)$ be a solution of (1.6) and suppose that $\|u(t)\|_{H^1} \leq R$ for $t \in [0,T]$. Then there exists a sequence $(\alpha_k)_{k\in\mathbb{N}}$, $\alpha_k \to 0$ as $k \to \infty$, such that*

$$\int_{|x|\geq k} |u(t,x)|^2 \,\mathrm{d}x \leq R^2 e^{-2\nu t} + \alpha_k \quad \text{for } t \in [0,T] \text{ and } k \in \mathbb{N}.$$

*Moreover, $\alpha_k$ depends only on $R$, $C$, $\nu$, $b(\cdot)$ and $c(\cdot)$.*

*Proof.* Let $\theta\colon \mathbb{R}_+ \to \mathbb{R}$ be a smooth function such that $0 \leq \theta(s) \leq 1$ for $s \in \mathbb{R}_+$, $\theta(s) = 0$ for $0 \leq s \leq 1$ and $\theta(s) = 1$ for $s \geq 2$. Let $D := \sup_{s\in\mathbb{R}_+} |\theta'(s)|$. Define $\theta_k(x) := \theta(|x|^2/k^2)$. Then, for $t \in [0,T]$, we have

$$\frac{d}{dt}\frac{1}{2}\int_{\mathbb{R}^n} \theta_k(x)|u(t,x)|^2 \,\mathrm{d}x = \int_{\mathbb{R}^n} \theta_k(x)u(t,x)u_t(t,x)\,\mathrm{d}x$$

$$= -\int_{\mathbb{R}^n} \nabla_x(\theta_k(x)u(t,x)) \cdot \nabla_x u(t,x)\,\mathrm{d}x + \int_{\mathbb{R}^n} \theta_k(x)u(t,x)F(x,u(t,x))\,\mathrm{d}x$$

Now we have

$$-\int_{\mathbb{R}^n} \nabla_x(\theta_k(x)u(t,x)) \cdot \nabla_x u(t,x)\,\mathrm{d}x$$

$$= -\int_{\mathbb{R}^n} \theta_k(x)|\nabla_x u(t,x)|^2 \,\mathrm{d}x - \frac{2}{k^2}\int_{\mathbb{R}^n} \theta'(|x|^2/k^2)u(t,x)\,x\cdot\nabla_x u(t,x)\,\mathrm{d}x$$

$$\leq \frac{2D}{k^2}\int_{k\leq|x|\leq\sqrt{2}k} |x|\,|u(t,x)|\,|\nabla_x u(t,x)|\,\mathrm{d}x \leq \frac{2\sqrt{2}D}{k}R^2.$$

On the other hand, by (2.1), by the Sobolev embedding $H^1 \hookrightarrow L^{2n/(n-2)}$ and by Hölder inequality, we have

$$\int_{\mathbb{R}^n} \theta_k(x)u(t,x)F(x,u(t,x))\,\mathrm{d}x$$

$$\leq -\nu\int_{\mathbb{R}^n} \theta_k(x)|u(t,x)|^2\,\mathrm{d}x + \int_{\mathbb{R}^n} \theta_k(x)b(x)|u(t,x)|^q\,\mathrm{d}x + \int_{\mathbb{R}^n} \theta_k(x)c(x)\,\mathrm{d}x$$

$$\leq -\nu\int_{\mathbb{R}^n} \theta_k(x)|u(t,x)|^2\mathrm{d}x + \left[\frac{(n-1)R}{(n-2)/2}\right]^q \left(\int_{|x|\geq k} |b(x)|^p\mathrm{d}x\right)^{1/p} + \int_{|x|\geq k}|c(x)|\mathrm{d}x.$$

Summing up, we have found a sequence $(\alpha_k)_{k\in\mathbb{N}}$, $\alpha_k \to 0$ as $k \to \infty$, such that

$$\frac{d}{dt}\int_{\mathbb{R}^n} \theta_k(x)|u(t,x)|^2 \,\mathrm{d}x \leq -2\nu\int_{\mathbb{R}^n} \theta_k(x)|u(t,x)|^2 \,\mathrm{d}x + \alpha_k.$$

Multiplying by $e^{2\nu t}$ and integrating on $[0,\bar{t}]$, we get

$$\int_{\mathbb{R}^n} \theta_k(x)|u(t,x)|^2 \,\mathrm{d}x \leq e^{-2\nu t}\int_{\mathbb{R}^n} \theta_k(x)|u(0,x)|^2 \,\mathrm{d}x + \alpha_k\frac{1}{2\nu}(1 - e^{-2\nu t}),$$

which in turn implies the thesis.  $\square$

The following proposition will allow us to recover $H^1$-admissibility from $L^2$-admissibility.

**Proposition 2.3.** *Let $G_j(x, u)$, $j = 1, 2, 3, \ldots$, and $G(x, u)$ be functions satisfying (1.3) with the same constant $C$. Assume that $G_j(x, u) \to G(x, u)$ for all $(x, u) \in \mathbb{R}^n \times \mathbb{R}$ and that $G_j(\cdot, 0) \to G(\cdot, 0)$ in $L^2(\mathbb{R}^n)$ as $j \to \infty$. Let $u_j \colon [0, T] \to H^1(\mathbb{R}^n)$ be a solution of (1.1) with $F := G_j$, $j = 1, 2, 3, \ldots$, and let $u \colon [0, T] \to H^1(\mathbb{R}^n)$ be a solution of (1.1) with $F := G$. Assume moreover that $\|u(t)\|_{H^1} \le R$ and $\|u_j(t)\|_{H^1} \le R$ for all $t \in [0, T]$ and all $j = 1, 2, 3, \ldots$. Finally, suppose that $u_j(0) \to u(0)$ in $L^2(\mathbb{R}^n)$ as $j \to \infty$. Then, for every $0 < \delta < T$, $u_j(t) \to u(t)$ in $H^1(\mathbb{R}^n)$ as $j \to \infty$, uniformly for $t \in [\delta, T]$. If $u_j(0) \to u(0)$ in $H^1(\mathbb{R}^n)$ as $j \to \infty$, then $u_j(t) \to u(t)$ in $H^1(\mathbb{R}^n)$ as $j \to \infty$, uniformly for $t \in [0, T]$.*

*Proof.* First, we observe that, if $u \in H^1(\mathbb{R}^n)$, then

$$G_j(x, u(x)) \to G(x, u(x)) \quad \text{as } j \to \infty$$

almost everywhere in $\mathbb{R}^n$. On the other hand, by (1.3),

$$|G_j(x, u(x)) - G(x, u(x))|^2$$
$$\le |G_j(x, 0) - G(x, 0)|^2 + |G_j(x, u(x)) - G_j(x, 0)|^2 + |G(x, u(x)) - G(x, 0)|^2$$
$$\le |G_j(x, 0) - G(x, 0)|^2 + 2C(|u(x)| + |u(x)|^{\beta+1})^2.$$

By the Dominated Convergence Theorem, we deduce that

$$(2.2) \qquad \hat{G}_j(u) \to \hat{G}(u) \quad \text{in } L^2(\mathbb{R}^n) \text{ as } j \to \infty, \text{ for any } u \in H^1(\mathbb{R}^n).$$

Assume now that $u_j(0) \to u(0)$ in $L^2(\mathbb{R}^n)$ as $j \to \infty$ and let $(t_j)_{j \in \mathbb{N}}$ be a sequence in $]0, T]$, converging to some $\bar{t} \in ]0, T]$. We procede like in the proof of Theorem 5.1 in [9]. Since

$$\|u_j(t_j) - u(\bar{t})\|_{H^1} \le \|u_j(t_j) - u(t_j)\|_{H^1} + \|u(t_j) - u(\bar{t})\|_{H^1}$$

and $\|u(t_j) - u(\bar{t})\|_{H^1} \to 0$ as $j \to \infty$, we need only to estimate $\|u_j(t_j) - u(t_j)\|_{H^1}$. Let $t \in ]0, T]$: then, in view of (1.7), we have

$$u_j(t) - u(t) = e^{-At}[u_j(0) - u(0)]$$
$$+ \int_0^t e^{-A(t-s)}[\hat{G}_j(u(s)) - \hat{G}(u(s))] \, ds + \int_0^t e^{-A(t-s)}[\hat{G}_j(u_j(s)) - \hat{G}(u_j(s))] \, ds$$

By (1.2), (1.4) and (1.5), we can find positive constants $K_1, K_2, \ldots$ (depending only on $\tilde{C}$, $M$, $R$ and $T$) such that

$$\|u_j(t) - u(t)\|_{H^1}$$
$$\le K_1 t^{-1/2} \|u_j(0) - u(0)\|_{L^2} + K_2 \int_0^t (t - s)^{-1/2} \|\hat{G}_j(u(s)) - \hat{G}(u(s))\|_{L^2} \, ds$$
$$+ K_3 \int_0^t (t - s)^{-1/2} \|u_j(s) - u(s)\|_{H^1} \, ds.$$

Setting

$$\gamma_j(t) := K_1 t^{-1/2} \|u_j(0) - u(0)\|_{L^2}$$

$$+ K_2 \int_0^t (t-s)^{-1/2} \|\hat{G}_j(u(s)) - \hat{G}(u(s))\|_{L^2} \, ds \quad \text{for } t \in ]0, T]$$

and applying Lemma 7.1.1 of [6], we obtain

$$\|u_j(t) - u(t)\|_{H^1} \le \gamma_j(t) + K_4 \int_0^t (t-s)^{-1/2} \gamma_j(s) \, ds.$$

Now, setting $\sigma := (\bar{t}/t_j)s$ for $s \in [0, t_j]$, we get

$$\|u_j(t_j) - u(t_j)\|_{H^1} \le \gamma_j(t_j) + K_4(t_j/\bar{t})^{1/2} \int_0^{\bar{t}} (\bar{t} - \sigma)^{-1/2} \gamma_j((t_j/\bar{t})\sigma) \, d\sigma.$$

Notice that $\gamma_j(t) \le K_5 t^{-1/2}$ for $t \in [0, T]$. The conclusion then follows by the Dominated Convergence Theorem, provided we can prove the following claim:

whenever $(\tau_j)_{j \in \mathbb{N}}$ is a sequence in $]0, T]$, converging to some $\bar{\tau} \in ]0, T]$, then $\gamma_j(\tau_j) \to 0$ as $j \to \infty$.

To this end, for $t, s \in [0, T]$, let us define

$$\chi_j(t, s) := \begin{cases} 0 & \text{if } t \le s \\ (t-s)^{-1/2} \|\hat{G}_j(u(s)) - \hat{G}(u(s))\|_{L^2} & \text{if } t > s \end{cases}$$

Then

$$\gamma_j(\tau_j) = K_1 \tau_j^{-1/2} \|u_j(0) - u(0)\|_{L^2} + K_2 \int_0^T \chi_j(\tau_j, s) \, ds.$$

If $s > \bar{\tau}$, then $s > \tau_j$ for all sufficiently large $j$, so $\chi_j(\tau_j, s) = 0$ for all sufficiently large $j$. Moreover, if $s < \bar{\tau}$, then $s < \tau_j$ for all sufficiently large $j$ and then, by (2.2), we have that $\chi_j(\tau_j, s) \to 0$ as $j \to \infty$. On the other hand, one has

$$|\chi_j(t, s)| \le K_5 |t - s|^{-1/2} \quad \text{for } t, s \in [0, T],$$

whence

$$|\chi_j(\tau_j, s)| \le K_5 |\tau_j - s|^{-1/2} \quad \text{for } s \in [0, T].$$

Since $\tau_j \to \bar{\tau}$ as $j \to \infty$, the sequence of functions

$$(|\tau_j - \cdot|^{-1/2})_{j \in \mathbb{N}}$$

converges to $|\bar{\tau} - \cdot|^{-1/2}$ in $L^1(0, T)$. By the Dominated Convergence Theorem, we have that $\gamma_j(\tau_j) \to 0$ as $j \to \infty$ and the claim is proved.

Finally, in order to complete the proof, we assume that $u_j(0) \to u(0)$ in $H^1(\mathbb{R}^n)$ as $j \to \infty$. In this case, for $t \geq 0$, we have

$$u_j(t) - u(0) = e^{-At} u_j(0) - u(0) + \int_0^t e^{-A(t-s)} \hat{G}(u_j(s)) \, \mathrm{d}s.$$

If $(t_j)_{j \in \mathbb{N}}$ is a sequence of positive numbers, $t_j \to 0$ as $j \to \infty$, then

$$\|u_j(t) - u(0)\|_{H^1} \leq \|e^{-At_j} u_j(0) - u(0)\|_{H^1} + K_6 \int_0^{t_j} (t_j - s)^{-1/2} \, \mathrm{d}s \to 0 \quad \text{as } j \to \infty,$$

and the proof is complete. $\square$

Now we are able to prove that all bounded subsets of $H^1(\mathbb{R}^n)$ are $\pi$-admissible. More precisely, we shall prove the following

**Theorem 2.4.** *Let $G_j(x, u)$, $j = 1, 2, 3, \dots$, and $G(x, u)$ be as in Proposition 2.3. Assume moreover that, for all $j = 1, 2, 3, \dots$, $G_j$ (and so also $G$) satisfies (2.1) (with $\nu$, $q$, $b$ and $c$ independent of $j$). Let $\pi_j$ (resp. $\pi$) be the local semiflow generated by equation (1.1) in $H^1(\mathbb{R}^n)$ with $F := G_j$ (resp. $F := G$). Finally, let $N$ be a closed bounded subset of $H^1(\mathbb{R}^n)$. Then $N$ is $\{\pi_j\}$-admissible.*

*Proof.* First, we choose $R > 0$ such that

$$N \subset B_{H^1}(R; 0) := \{ u \in H^1(\mathbb{R}^n) \mid \|u\|_{H^1} \leq R \}.$$

Now let $(u_j)_{j \in \mathbb{N}}$ be a sequence in $H^1(\mathbb{R}^n)$ and let $(t_j)_{j \in \mathbb{N}}$ be a sequence of positive numbers such that $t_j \to \infty$ as $j \to \infty$ and $u_j \pi_j[0, t_j] \subset N$ for all $j \in \mathbb{N}$.

By carefully checking the proof of Theorem 3.3.3 in [6], we find out that there exists $\tau > 0$ such that, for all $u \in B_{H^1}(R; 0)$, $u\pi t$ is defined for $t \in [0, \tau]$ and $u\pi[0, \tau] \subset B_{H^1}(2R; 0)$.

Without loss of generality, we can assume that $t_j > \tau$ for all $j \in \mathbb{N}$. Since

$$\|u_j \pi_j(t_j - \tau)\|_{H^1} \leq R \quad j = 1, 2, 3, \dots,$$

then there exists $v \in H^1(\mathbb{R}^n)$, $\|v\|_{H^1} \leq R$, such that, up to a subsequence,

$$(2.3) \qquad\qquad u_j \pi_j(t_j - \tau) \rightharpoonup v \quad \text{in } H^1(\mathbb{R}^n).$$

Notice that $u_j \pi_j t$ is defined for $t \in [0, \tau]$ and $u_j \pi_j[0, \tau] \subset B_{H^1}(R; 0)$. Besides, $v\pi t$ is defined for $t \in [0, \tau]$ and $v\pi[0, \tau] \subset B_{H^1}(2R; 0)$.

Let $k \in \mathbb{N}$ and $\theta_k$ be as in the proof of Proposition 2.2. By Proposition 2.2 we have

$$\int_{\mathbb{R}^n} \theta_k(x) |u_j \pi_j(t_j - \tau)(x)|^2 \, \mathrm{d}x \leq R^2 e^{-2\nu(t_j - \tau)} + \alpha_k.$$

Let $\epsilon > 0$ be fixed. Take $k$ and $j_0$ so large that, for all $j \geq j_0$, one has $R^2 e^{-2\nu(t_j - \tau)} + \alpha_k \leq \epsilon$. Then

$$\{ u_j \pi_j(t_j - \tau) \mid j \geq j_0 \} = \{ \theta_k [u_j \pi_j(t_j - \tau)] + (1 - \theta_k) [u_j \pi_j(t_j - \tau)] \mid j \geq j_0 \}$$
$$\subset \{ \theta_k [u_j \pi_j(t_j - \tau)] \mid j \geq j_0 \} + \{ (1 - \theta_k) [u_j \pi_j(t_j - \tau)] \mid j \geq j_0 \}$$
$$(2.4) \qquad \subset B_{L^2}(\epsilon; 0) + \{ (1 - \theta_k) [u_j \pi_j(t_j - \tau)] \mid j \geq j_0 \}.$$

The set

$$\{ (1 - \theta_k) [u_j \pi_j(t_j - \tau)] \mid j \geq j_0 \}$$

consists of functions of $H^1(\mathbb{R}^n)$ which are equal to zero outside the ball $B_{\sqrt{2}k}(0)$ in $\mathbb{R}^n$. On the other hand, the $H^1$ norm of these functions is bounded by a constant depending only on $R$ and $D$. Then, by Rellich Theorem, we deduce that this set is precompact in $L^2(\mathbb{R}^n)$. Hence we can cover it by a finite number of balls of radius $\epsilon$ in $L^2(\mathbb{R}^n)$. This observation, together with (2.4), implies that the set $\{ u_j \pi_j(t_j - \tau) \mid j \geq j_0 \}$ is totally bounded and hence precompact in $L^2(\mathbb{R}^n)$. Thus, up to a subsequence, we can assume that

$$(2.5) \qquad u_j \pi_j(t_j - \tau) \to v \quad \text{in } L^2(\mathbb{R}^n).$$

Finally, by Proposition 2.3, we have that, up to a subsequence,

$$(2.6) \qquad u_j \pi_j t_j = [u_j \pi_j(t_j - \tau)] \pi_j \tau \to v \pi \tau \quad \text{in } H^1(\mathbb{R}^n)$$

as $j \to \infty$, and the theorem is proved. $\square$

We end this section by stating and proving an important consequence of Theorem 2.4. First, we give the following

**Definition 2.5.** *We say that a function $\sigma \colon \mathbb{R} \to H^1(\mathbb{R}^n)$ is a full solution of equation (1.1) iff $\sigma(t) = \sigma(s)\pi(t - s)$ for all $t \geq s$, where $\pi$ is the local semiflow generated by (1.1) in $H^1(\mathbb{R}^n)$.*

Now we have:

**Corollary 2.6.** *Let $G_j(x, u)$, $j = 1, 2, 3, \ldots$, and $G(x, u)$ be as in Theorem 2.4 and let $\pi_j$ (resp. $\pi$) be the local semiflow generated by equation (1.1) in $H^1(\mathbb{R}^n)$ with $F := G_j$ (resp. $F := G$). Let $R$ be a positive constant and, for all $j \in \mathbb{N}$, let $\sigma_j \colon \mathbb{R} \to H^1(\mathbb{R}^n)$ be a full solution of (1.1) with $F := G_j$ such that*

$$\sup_{t \in \mathbb{R}} \|\sigma_j(t)\|_{H^1} \leq R.$$

*Under these hypotheses, there exists a subsequence of $(\sigma_j)_{j \in \mathbb{N}}$, again denoted by $(\sigma_j)_{j \in \mathbb{N}}$, and a full solution $\sigma \colon \mathbb{R} \to H^1(\mathbb{R}^n)$ of (1.1) with $F := G$, such that*

$$\sigma_j(t) \to \sigma(t) \quad \text{as } j \to \infty$$

*uniformly on every bounded subinterval of* $\mathbb{R}$.

*Proof.* As in the proof of Theorem 2.4, we begin by taking $\tau > 0$ such that, for all $u \in B_{H^1}(R; 0)$, $u\pi t$ is defined for $t \in [0, \tau]$ and $u\pi[0, \tau] \subset B_{H^1}(2R; 0)$. Then we fix once and for all a sequence $(t_j)_{j\in\mathbb{N}}$ of positive numbers, with $t_j \to \infty$ as $j \to \infty$.

Let $k \in \mathbb{Z}$. For all sufficiently large $j$, we have

$$\sigma_j(k\tau) = \sigma_j(k\tau - t_j)\pi_j t_j.$$

Then, by Theorem 2.4, there is a subsequence of $(\sigma_j(k\tau))_{j\in\mathbb{N}}$, again denoted by $(\sigma_j(k\tau))_{j\in\mathbb{N}}$, and there exists $v(k\tau) \in H^1(\mathbb{R}^n)$ such that $\sigma_j(k\tau)$ converges strongly to $v(k\tau)$ in $H^1(\mathbb{R}^n)$ as $j \to \infty$. In particular, $\|v(k\tau)\|_{H^1} \leq R$. Using Cantor's diagonal procedure we obtain the existence of a subsequence of $(\sigma_j)_{j\in\mathbb{N}}$, again denoted by $(\sigma_j)_{j\in\mathbb{N}}$, and a sequence $v(k\tau) \in H^1(\mathbb{R}^n)$, $k \in \mathbb{Z}$, such that, for every $k \in \mathbb{Z}$,

$$\sigma_j(k\tau) \to v(k\tau) \quad \text{in } H^1(\mathbb{R}^n) \text{ as } j \to \infty.$$

By Proposition 2.3, we have that, for all $k \in \mathbb{Z}$,

$$\sigma_j(k\tau)\pi_j t \to v(k\tau)\pi t \quad \text{in } H^1(\mathbb{R}^n) \text{ as } j \to \infty, \text{ uniformly on } [0, \tau].$$

In particular, one has $\sigma_j(k\tau)\pi_j\tau \to v(k\tau)\pi\tau$. On the other hand, $\sigma_j(k\tau)\pi_j\tau = \sigma_j((k+1)\tau) \to v((k+1)\tau)$. Hence we deduce that $v((k+1)\tau) = v(k\tau)\pi\tau$ for all $k \in \mathbb{Z}$. We can therefore define

$$\sigma(t) := v(k\tau)\pi(t - k\tau) \quad \text{for } t \in [k\tau, (k+1)\tau],$$

which is easily seen to be a full solution of of (1.1) with $F := G$. Moreover,

$$\sigma_j(t) \to \sigma(t) \quad \text{as } j \to \infty$$

uniformly on every bounded subinterval of $\mathbb{R}$.   $\square$

## 3. ASYMPTOTICALLY LINEAR EQUATIONS

In this section we concentrate on equations which behave linearly at infinity. More precisely, we assume that $F(x, u)$ satisfies (1.3) with $\beta = 0$ and (2.1) with $q = 2$. Moreover, we assume that

$$(3.1) \qquad \lim_{|u|\to\infty} \frac{F(x, u)}{u} = \alpha(x) := -\alpha_1(x) + \alpha_2(x) \quad \text{for all } x \in \mathbb{R}^n,$$

where $\alpha_1 \in L^\infty(\mathbb{R}^n)$, with $\alpha_1(x) \geq \tilde{\nu} > 0$ for all $x \in \mathbb{R}^n$, and $\alpha_2 \in L^\rho(\mathbb{R}^n)$, with $n \leq \rho < \infty$.

*Remark.* The results of this section still hold if in (3.1) (resp. in (3.3)) one makes the alternative assumption that $\alpha_2 \in L^\infty(\mathbb{R}^n)$ and $\lim_{k\to\infty} \text{ess sup}_{|x|\geq k} |\alpha_2(x)| = 0$ (resp. $\gamma_2 \in L^\infty(\mathbb{R}^n)$ and $\lim_{k\to\infty} \text{ess sup}_{|x|\geq k} |\gamma_2(x)| = 0$). In fact, one only needs to slightly modify the proof of Lemma 3.1 below.

Our aim is to extend the results of [12] to equations on $\mathbb{R}^n$. First, observe that the spectrum (and hence *a fortiori* the essential spectrum) of the operator $-\Delta + \alpha_1(\cdot)$ is contained in the interval $[\tilde{\nu}, +\infty[$. On the other hand, the multiplication operator $u(\cdot) \mapsto \alpha_2(\cdot)u(\cdot)$, defined on $H^1(\mathbb{R}^n)$ with values in $L^2(\mathbb{R}^n)$, turns out to be relatively compact with respect to $-\Delta + \alpha_1(\cdot)$. In fact, one has the following

**Lemma 3.1.** *The operator*

$$\alpha_2 \circ (-\Delta + \alpha_1)^{-1} \colon L^2(\mathbb{R}^n) \to L^2(\mathbb{R}^n)$$

*is compact.*

*Proof.* Since $(-\Delta + \alpha_1)^{-1} \colon L^2(\mathbb{R}^n) \to H^1(\mathbb{R}^n)$ is bounded, it is enough to take a bounded sequence $(u_j)_{j\in\mathbb{N}}$ in $H^1(\mathbb{R}^n)$ and show that, up to a subsequence, $(\alpha_2 u_j)_{j\in\mathbb{N}}$ converges in $L^2(\mathbb{R}^n)$. For $k > 0$, one has

$$\int_{|x|\geq k} [\alpha_2(x)u_j(x)]^2 \mathrm{d}x \leq \left(\int_{|x|\geq k} \alpha_2(x)^\rho \mathrm{d}x\right)^{2/\rho} \left(\int_{|x|\geq k} u_j(x)^{2\rho/(\rho-2)} \mathrm{d}x\right)^{(\rho-2)/\rho}.$$

Since $2 \leq 2\rho/(\rho-2) \leq 2n/(n-2)$, given $\epsilon > 0$, we can choose $k$ so large that

$$\int_{|x|\geq k} [\alpha_2(x)u_j(x)]^2 \mathrm{d}x \leq \epsilon$$

for all $j \in \mathbb{N}$. Then we proceed like in the proof of Theorem 2.4. If $\chi_k$ is the characteristic function of $\{\, |x| \geq k \,\}$, we have

$$\{\, \alpha_2 u_j \mid j \in \mathbb{N} \,\} = \{\, \chi_k \alpha_2 u_j + (1-\chi_k)\alpha_2 u_j \mid j \in \mathbb{N} \,\}$$
$$\subset \{\, \chi_k \alpha_2 u_j \mid j \in \mathbb{N} \,\} + \{\, (1-\chi_k)\alpha_2 u_j \mid j \in \mathbb{N} \,\}$$
$$\subset B_{L^2}(\epsilon; 0) + \{\, (1-\chi_k)\alpha_2 u_j \mid j \in \mathbb{N} \,\}$$

The functions $u_j$, $j \in \mathbb{N}$, restricted to $\{\, |x| < k \,\}$, form a bounded subset of $H^1(\{\, |x| < k \,\})$. By Rellich's Theorem, up to a subsequence, they converge in $L^2(\{\, |x| < k \,\})$. Since $\alpha_2 \in L^\infty(\mathbb{R}^n)$, it follows that, up to a subsequence, also the functions $(1-\chi_k)\alpha_2 u_j$ converge in $L^2(\mathbb{R}^n)$. Hence, the set $\{\, (1-\chi_k)\alpha_2 u_j \mid j \in \mathbb{N} \,\}$ is precompact in $L^2(\mathbb{R}^n)$: so we can cover it by a finite number of balls of radius $\epsilon$ in $L^2(\mathbb{R}^n)$. It follows that the set $\{\, \alpha_2 u_j \mid j \in \mathbb{N} \,\}$ is totally bounded and hence precompact in $L^2(\mathbb{R}^n)$. $\square$

By Weyl's Theorem (see e.g. [10, p. 113]), also the essential spectrum of the operator $-\Delta + \alpha_1(\cdot) - \alpha_2(\cdot)$ is contained in $[\tilde{\nu}, +\infty[$. In particular, the part of the spectrum of $-\Delta + \alpha_1(\cdot) - \alpha_2(\cdot)$ contained in $]-\infty, \tilde{\nu}/2[$ is a finite set, consisting of isolated eigenvalues with finite multiplicity. From now on, we assume that the following *non-resonance condition at infinity* is satisfied:

$$(3.2) \qquad\qquad \ker(-\Delta + \alpha_1(\cdot) - \alpha_2(\cdot)) = (0).$$

Whenever $F$ satisfies (1.3), we denote by $\pi_F$ the semiflow generated by (1.1) and by $\mathcal{K}_F$ the union of all bounded full orbits of $\pi_F$. We have the following

**Proposition 3.2.** *Assume that F satisfies (1.3), (2.1), (3.1) and (3.2). Write* $B(x,u) := \alpha(x)u$. *Then there exists* $R > 0$ *such that, for all* $\lambda \in [0,1]$, $\mathcal{K}_{\lambda F + (1-\lambda)B}$ *is contained in* $B_{H^1}(R;0)$. *Moreover,* $\mathcal{K}_{\lambda F + (1-\lambda)B}$ *is compact in* $H^1(\mathbb{R}^n)$.

*Proof.* The proof is by contradiction. Suppose the theorem is not true; then there are a sequence $(\lambda_j)_{j \in \mathbb{N}}$ in $[0,1]$, a sequence $(\sigma_j)_{j \in \mathbb{N}}$ of bounded full solutions of $\pi_{\lambda_j F + (1-\lambda_j)B}$ and a sequence of positive numbers $(R_j)_{j \to \infty}$, $R_j \to \infty$ as $j \to \infty$, such that

$$R_j = \sup_{t \in \mathbb{R}} \|\sigma_j(t)\|_{H^1} \quad \text{and} \quad \|\sigma_j(0)\|_{H^1} \geq R_j - 1 > 0.$$

Let us define

$$G_j(x,u) := \lambda_j R_j^{-1} F(x, R_j u) + (1 - \lambda_j) R_j^{-1} B(x, R_j u)$$
$$= \alpha(x)u + \lambda_j \left( R_j^{-1} F(x, R_j u) - \alpha(x)u \right)$$

and

$$\zeta_j(t) := R_j^{-1} \sigma_j(t).$$

Notice that $\zeta_j$ is a bounded full solution of $\pi_{G_j}$, with $\sup_{t \in \mathbb{R}} \|\zeta_j(t)\|_{H^1} = 1$ and $\|\zeta_j(0)\|_{H^1} \to 1$ as $j \to \infty$.

It is easy to check that $G_j$ and $B$ satisfy (1.3) and (2.1) uniformly with respect to $j \in \mathbb{N}$. Moreover, by (3.1), $G_j(x,u) \to B(x,u)$ for all $(x,u) \in \mathbb{R}^n \times \mathbb{R}$ and $G_j(\cdot, 0) \to B(\cdot, 0)$ in $L^2(\mathbb{R}^n)$ as $j \to \infty$.

Then, by Corollary 2.6, there exists a subsequence of $(\zeta_j)_{j \in \mathbb{N}}$, again denoted by $(\zeta_j)_{j \in \mathbb{N}}$, and a full solution $\zeta : \mathbb{R} \to H^1(\mathbb{R}^n)$ of $\pi_B$, such that

$$\zeta_j(t) \to \zeta(t) \quad \text{as } j \to \infty$$

uniformly on every bounded subinterval of $\mathbb{R}$. In particular, $\|\zeta(0)\|_{H^1} = 1$. On the other hand, by [13, Th. I.11.1], the only bounded full solution of $\pi_B$ is $\zeta(t) \equiv 0$, a contradiction.

Finally, the fact that $\mathcal{K}_{\lambda F + (1-\lambda)B}$ is compact in $H^1(\mathbb{R}^n)$ is a straightforward consequence of Corollary 2.6.   $\square$

We denote by $m$ the total multiplicity of the negative eigenvalues of $-\Delta + \alpha_1(\cdot) - \alpha_2(\cdot)$. By [13, Th. I.11.1], we have $h(\pi_B, \{0\}) = \Sigma^m$. By propositions 2.3 and 3.2, by Theorem 2.4 and by the Continuation Theorem [13, Th. I.12.2], we finally obtain the following extension of Theorem 2.2 of [12]:

**Theorem 3.3.** *Assume that F satisfies (1.3), (2.1), (3.1) and (3.2). Let $m$ be the total multiplicity of the negative eigenvalues of* $-\Delta + \alpha_1(\cdot) - \alpha_2(\cdot)$. *Then*

$$h(\pi_F, \mathcal{K}_F) = \Sigma^m.$$

*In particular, $\mathcal{K}_F$ is non-empty and irreducible.*   $\square$

As a consequence, we have the following

**Corollary 3.4.** *Assume that $F$ satisfies (1.3), (2.1), (3.1) and (3.2). Then there exists at least one equilibrium of equation (1.1).*

*Proof.* This is simple consequence of the gradient-like structure of the semiflow $\pi_F$. In fact, if $P(x,u) := \int_0^u F(x,s)\,\mathrm{d}s$, then

$$V(u) := \frac{1}{2}\int_{\mathbb{R}^n} |\nabla u(x)|^2\,\mathrm{d}x - \int_{\mathbb{R}^n} P(x,u(x))\,\mathrm{d}x$$

is a $C^1$ Lyapunov functional for $\pi_F$, such that

$$\frac{d}{dt}V(u(t)) = -|\dot{u}(t)|^2$$

along any solution $u(t)$ of (1.1). It is a well known fact that then the $\alpha$- and $\omega$-limit of a bounded full solution of (1.1) are non-empty, compact, connected and consist only of equilibria of (1.1) (see e.g. [5], [7] or [16]). $\square$

*Remark.* The result of Corollary 3.4 is not new, as it was already proved in [15] under conditions which are essentially equivalent to (2.1), (3.1) and (3.2). Actually, the result of [15] is even more general than ours, since it does not require differentiability but only continuity of $F$ with respect to $u$. However, if $F(x,0) \equiv 0$ then 0 is an equilibrium of (1.1) and a result like that of Corollary 3.4 or [15, Th. 1] is not very informative. If this is the case, our dynamical approach allows us to obtain more precise existence results for nontrivial equilibria of (1.1), by analysing the linearization of (1.1) at 0, as we shall explain below. Moreover, we can even prove the existence of heteroclinic connections between 0 and some set of nontrivial equilibria.

In the final part of this section, we assume that $F(x,0) \equiv 0$, so 0 is an equilibrium of (1.1). We assume that $F(x,u)$ satisfies (1.3) with $\beta = 0$ and (2.1) with $q = 2$. Setting $F'_u(x,0) =: \gamma(x)$, we assume that

$$(3.3) \qquad\qquad \gamma(x) := -\gamma_1(x) + \gamma_2(x)$$

where $\gamma_1 \in L^\infty(\mathbb{R}^n)$, with $\gamma_1(x) \geq \tilde{\nu} > 0$ for all $x \in \mathbb{R}^n$, and $\gamma_2 \in L^\rho(\mathbb{R}^n)$, with $n \leq \rho < \infty$. Again, the spectrum (and hence *a fortiori* the essential spectrum) of the operator $-\Delta + \gamma_1(\cdot)$ is contained in the interval $[\tilde{\nu}, +\infty[$. On the other hand, the multiplication operator $u(\cdot) \mapsto \gamma_2(\cdot)u(\cdot)$, defined on $H^1(\mathbb{R}^n)$ with values in $L^2(\mathbb{R}^n)$, is relatively compact with respect to $-\Delta + \gamma_1(\cdot)$. By Weyl's Theorem, also the essential spectrum of the operator $-\Delta + \gamma_1(\cdot) - \gamma_2(\cdot)$ is contained in $[\tilde{\nu}, +\infty[$. In particular, the part of the spectrum of $-\Delta + \gamma_1(\cdot) - \gamma_2(\cdot)$ contained in $]-\infty, \tilde{\nu}/2[$ is a finite set, consisting of isolated eigenvalues with finite multiplicity. We assume that the following *non-resonance condition at zero* is satisfied:

$$(3.4) \qquad\qquad \ker(-\Delta + \gamma_1(\cdot) - \gamma_2(\cdot)) = (0).$$

Then we have the following result, which is 'dual' to Theorem 3.2:

**Proposition 3.5.** *Assume that $F$ satisfies (1.3), (2.1), (3.3) and (3.4). Write $D(x, u) := \gamma(x)u$. Then there exists $\delta > 0$ such that, for all $\lambda \in [0, 1]$, $\{0\}$ is the maximal invariant set of $\pi_{\lambda F + (1-\lambda)D}$ in $B_{H^1}(\delta; 0)$.*

*Proof.* The proof is by contradiction and is completely analogous to that of Theorem 3.2; therefore it is left to the reader. $\quad\square$

We denote by $m'$ the total multiplicity of the negative eigenvalues of $-\Delta + \gamma_1(\cdot) - \gamma_2(\cdot)$. By [13, Th. I.11.1], we have $h(\pi_D, \{0\}) = \Sigma^{m'}$. By propositions 2.3 and 3.5, by Theorem 2.4 and by the Continuation Theorem [13, Th. I.12.2], we finally obtain the following extension of Theorem 3.1 of [12]:

**Theorem 3.6.** *Assume that $F$ satisfies (1.3), (2.1), (3.3) and (3.4). Let $m'$ be the total multiplicity of the negative eigenvalues of $-\Delta + \gamma_1(\cdot) - \gamma_2(\cdot)$. Then*

$$h(\pi_F, \{0\}) = \Sigma^{m'}.$$

$\square$

Finally, we have the following extension of Theorem 3.3 of [12]:

**Corollary 3.7.** *Assume that $F$ satisfies (1.3), (2.1), (3.1) and (3.2). Besides, assume that $F(x, 0) \equiv 0$ and $F$ satisfies (3.3) and (3.4). Let $m$ be the total multiplicity of the negative eigenvalues of $-\Delta + \alpha_1(\cdot) - \alpha_2(\cdot)$ and let $m'$ be the total multiplicity of the negative eigenvalues of $-\Delta + \gamma_1(\cdot) - \gamma_2(\cdot)$. If $m \neq m'$, there exists at least one nontrivial equilibrium of equation (1.1). Moreover, there exists a heteroclinic orbit of $\pi_F$ connecting $0$ with a set of nontrivial equilibria.*

*Proof.* The proof is analogous to that of Theorem 3.3 of [12], to which the reader is referred. $\quad\square$

Martino Prizzi, Università degli Studi di Trieste, Dipartimento di Scienze Matematiche, via Valerio 12, 34127 Trieste, Italy

*E-mail address*: `prizzi@mathsun1.univ.trieste.it`